\documentclass[default]{sn-jnl}

\jyear{2022}%

\usepackage{graphicx}
\usepackage{float}
\usepackage{amsmath,amsxtra,amssymb,latexsym, amscd,amsthm}
\usepackage{makecell}

\usepackage{natbib}  
\setcitestyle{authoryear,open={(},close={)}}

\usepackage{caption}
\usepackage{algorithm}

\DeclareCaptionFormat{myformat}{#3}
\DeclareMathOperator\dom{dom}
\DeclareMathOperator\argmin{arg\,min}

\theoremstyle{thmstyleone}%
\newtheorem{theorem}{Theorem}
\newtheorem{proposition}[theorem]{Proposition}%

\theoremstyle{thmstyletwo}%
\newtheorem{example}{Example}%
\newtheorem{remark}{Remark}%

\theoremstyle{thmstylethree}%

\begin{document}

\title[Solving the problem of batch deletion and insertion members in LKH by DCA]{Solving the problem of batch deletion and insertion members in the Logical Key Hierarchy structure by a DC Programming approach}


\author[1,2]{\fnm{Hoai An} \sur{Le Thi}}\email{hoai-an.le-thi@univ-lorraine.fr}

\author*[1]{\fnm{Thi Tuyet Trinh} \sur{Nguyen}}\email{thi-tuyet-trinh.nguyen@univ-lorraine.fr}

\affil*[1]{\orgdiv{Universit\'e de Lorraine}, \orgname{LGIPM, D\'epartement IA}, \orgaddress{\postcode{F-57000}, \city{Metz}, \country{France}}}

\affil[2]{\orgname{Institut Universitaire de France (IUF)}}


\abstract{In secure group communications, users of a group share a common group key to prevent eavesdropping and protect the exchange content. A key server distributes the group key as well as performs group rekeying whenever the membership changes dynamically. Instead of rekeying after each join or leave request, we use batch rekeying to alleviate the out-of-sync problem and improve the efficiency. In this paper, we propose an optimization approach to the problem of updating group key in the Logical Key Hierarchy (LKH) structure with batch rekeying. A subtree of new nodes can be appended below a leaf node or is replaced the position of leaving node on the binary key tree. The latter has a lower updating key cost than the former since when a member leaves, all the keys on the path from the root to the deletion node must be updated anyway. We aim to minimize the total rekeying cost, which is the cost of deletion and insertion members while keeping the tree as balanced as possible. The mentioned problem is represented by a unified (deterministic) optimization model whose objective function contains discontinuous step functions with binary variables. Thanks to an exact penalty technique, the problem is equivalently reformulated as a standard DC (Difference of Convex functions) program that can be solved efficiently by DCA (DC algorithm). Numerical experiments have been studied intensively to justify the merit of our proposed approach as well as the corresponding DCA.}




\maketitle

\section{Introduction}\label{sec1}

Many secure group communication systems are based on a group key that is secretly shared by group members. In order to provide security for such communications, the existing systems encrypt the data using this group key and send the corresponding ciphertext to all members. Therefore, securing group communications (i.e., providing confidentiality, authenticity, and integrity of messages delivered between group members) has become an important Internet design issue. Due to the dynamic nature of group membership, it becomes necessary to change the group key in an efficient and secure manner when members join or leave the group.

A group key management system must specifically assure two forms of secrecy: backward and forward secrecy. Backward secrecy means that if the system changes the group key after a new user joins, the new member will be unable to decrypt previous group messages. The departing user will not be able to access future group communication if the system rekeys after a current user departs or is ejected from the system; this is known as forward secrecy. Key server can work in three different batch modes \citep{yang2001reliable} to meet the requirements of different applications and find a balance between performance and group security. The first method is batch rekeying, in which the key server handles both join and leave requests periodically in a batch. The second method is batch leave rekeying, in which the key server processes each join request instantly to shorten the time it takes for a new user to access group communications, but batch leave requests. The last mode is batch join rekeying, in which the key server processes each leave request quickly to avoid exposure to departed users, but processes join requests in batches.

In \cite{le2022dc}, we proposed an approach to an important problem in centralized dynamic group key management that consists in finding a set of leaf nodes in a binary key tree to insert new members while minimizing the insertion cost and keeping the tree as balanced as possible. We opt for the batch inserting technique and individual deletion. In actuality, processing each leave request immediately ensure the forward secrecy for group communications. However, the number of updated keys generated by individual leave rekeying can be much more than the sum of those generated by batch rekeying. In this case, key server and group members are assigned to perform more rekeying operations rather than making use of the service. Batch rekeying \citep{li2001batch,zhang2003protocol,ng2007dynamic,pegueroles2003balanced,vijayakumar2012rotation} can alleviate the out-of-sync problem and improve efficiency. \cite{li2001batch} suggested a method for updating the key tree and generating a rekey subtree containing a collection of join and depart requests at the end of each rekey period. \cite{ng2007dynamic} developed the batch balanced approach, which accounts for both new and leaving members in order to reduce key tree height differences. \cite{vijayakumar2012rotation} introduced the use of rotation-based key tree methods to maintain the tree's balance even when batch depart requests outnumber batch join requests. Existing techniques do not account for both key tree balance and rekeying costs concurrently, resulting in an unbalanced key tree or high rekeying costs.

When the membership of a group changes dynamically, the group key and all related keys on the binary tree must be updated. In batch rekeying, the key server aggregates the total number of joining and leaving members over a specified time interval and then changes the associated keys. Therefore, batch rekeying techniques increase efficiency in the number of required messages and it takes advantage of the possible overlap of new keys for multiple rekeying requests, and then reduces the possibility of generating redundant new keys. In actuality, the deletion nodes are also the leaf nodes in the tree and a subtree of new nodes can be appended below a leaf node or is replaced the position of leaving node on the binary key tree. The latter has a lower updating key cost than the former since when a member leaves, all the keys on the path from the root to the deletion node must be updated anyway. This is the advantage of this model compared to the existing algorithms.  
 
\textbf{Our contributions.} We develop an optimization approach to the problem of updating group key in the LKH structure with batch rekeying technique. As the deleting and inserting nodes should be leaf nodes, the first step of our algorithm consists in finding the set of all leaf nodes of the given tree. In the second stage, based on the found leaf nodes, we find a set of leaf nodes to delete leaving members and insert new members while minimizing the total rekeying cost, which is the cost of deletion and insertion members, and taking the balance of the tree into account. The proposed optimization model is a combinatorial problem with discontinuous step functions and binary variables. It is very challenging to handle such kinds of programs by standard methods where the source of difficulty comes from the discontinuity of the objective and the binary nature of the solutions. The problem is first equivalently formulated to remove the discontinuity of the objective, and then the latter problem is reformulated as a DC program via an exact penalty technique, where the DCA is at our disposal as an efficient algorithm in DC programming \citep{dinh2010efficient}.

The rest of the paper is organized as follows. Section 2 gives an optimization approach to the problem of updating group key in the LKH structure with batch rekeying. Section 3 presents the solution method based on DC programming and DCA. The implementation of the algorithm for solving the problem in LKH structure and numerical experiments are presented in Section 4. Finally, some conclusions are provided in Section 5.

\section{An optimization approach to the problem of batch deletion and insertion members}\label{sec2}

\subsection{Problem description}\label{subsec2.1}

We consider the problem of batch deletion and insertion members in the LKH structure. When the membership of a group changes dynamically, the group key and all related keys on the binary tree must be updated. Our main goal here is to find a set of leaf nodes in a binary key tree to delete leaving members and insert new members while minimizing the rekeying cost (the number of updated keys) and keeping the tree as balanced as possible.

\subsection{A two-step algorithm and the optimization model}\label{subsec2.2}

As the deletion and insertion nodes should be leaf nodes, we first find all the leaf nodes in a full binary tree. Second, we propose an optimization model based on the found leaf nodes that minimizes the cost of deleting departure members and adding new members while maintaining the tree's balance. The following describes the optimization model.

Given a full binary tree that can be represented as an ordered set $T=\{2^a+b$: there exists a node at the $b$-th horizontal position of level $a$ on the tree, $0\leq a\leq h,0\leq b\leq 2^a-1,a,b\in \mathbb{N}\}$, where $h$ is the height of tree, we seek a set of the leaf nodes $L\subseteq T$ that is defined as $L=\{t\in T: 2t\notin T,2t+1\notin T\},l=\vert L\vert$. $A$ is a set of leaving members and $A \subseteq L,l_2=\vert A\vert$. The set of remaining leaf nodes on the tree $T$ is denoted as $L'=L\setminus A, l_1=\vert L\setminus A\vert$. For convenience, we sort $T,L,A,L'$ in ascending order. The distance from the root to the leaf node $L'[i]$ and the departure leaf node $A[k]$ are given by 
$d_i = \lfloor\log_2 L'[i]\rfloor,i=\overline{1,l_1}$ and $d_k =\lfloor\log_2 A[k]\rfloor,k=\overline{1,l_2}$, respectively. Let the set of joining members be $M=\{1,2,…,m\}$.

\begin{figure}[h]%
\centering
\includegraphics[width=0.7\textwidth]{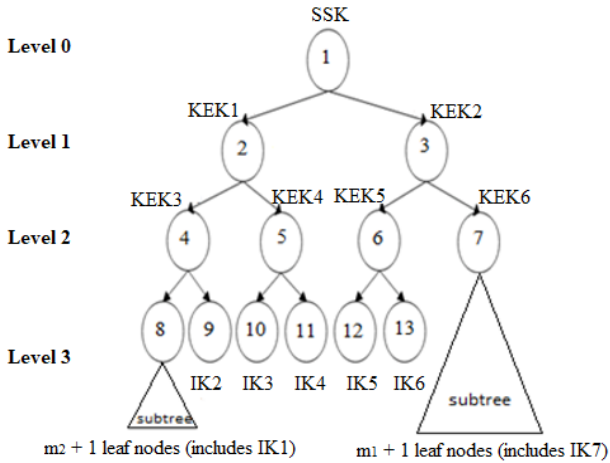}
\caption{The tree structure after deleting departure nodes and inserting new nodes.}\label{fig1}
\end{figure}

\begin{remark}

i) Each leaving node will be either replaced with a certain number of new members or only deleted from the key tree (illustrated in Fig. \ref{fig1}).  If a subtree of new nodes is replaced the position of leaving nodes, the rekeying cost gets reduced by the number of all keys on the path from the root to the deletion node which are certainly required to be updated when a member assigned to this node is leaving from the group. 

ii) After inserting $m$ members into the tree, each leaf node will be either assigned with a certain number of new members or not changed (illustrated in Fig. \ref{fig1}). If a leaf node $L[i]$ is assigned some new members, it becomes a key interior node, otherwise, $L[i]$ remains a leaf node.

iii) A valid insertion procedure will create a full binary subtree below each assigned leaf node. According to Handshaking lemma \cite{vasudev2006graph} and the proprieties of a full binary tree, we have the following assertion: The cost of appending new nodes at the position of leaf node on the full binary tree only depends on the number of new nodes to be inserted at that leaf node but it does not depend on the configuration of the full binary subtree.
\end{remark}

Now we propose the optimization model for minimizing the key updating cost. Let $x_{ij}$ be binary variables defined by $x_{ij} = 1$ if a new member $j \in M$ is inserted into the subtree below the leaf node $L'[i]$,  $x_{ij} = 0$ otherwise, $i=\overline{1,l_1},j=\overline{1,m}$. Let $y_{kj}$ be binary variables defined by $y_{kj}=1$ if a leaf node $A[k]$ is deleted and a new member $j \in M$ is inserted into the subtree at the leaf node $A[k]$,  $y_{kj} = 0$ if the leaf node $A[k]$ is only deleted, $k=\overline{1,l_2},j=\overline{1,m}$.

Since every new member is appended below only one leaf node of the original tree, $\sum_{i=1}^{l_1}{x_{ij}}+\sum_{k=1}^{l_2}{y_{kj}} = 1$, $\forall j = \overline{1,m}$. For a given value $i\in \{1,2,...,l_1\}$, the number of new nodes $j$ being inserted at the leaf node $L'[i]$ is calculated as $m_i= \sum_{j=1}^{m} x_{ij}$. It means that a subtree at the leaf node $L'[i]$ includes $m_i+1$ leaf nodes (including $m_i$ new nodes and the old leaf node $L'[i]$) as illustrated in Fig. \ref{fig1}. For the leaf node $L'[i]$ that is chosen to append some new members (equivalent to $\sum_{j=1}^m x_{ij}>0$), the cost to build a corresponding subtree is $2\sum_{j=1}^m x_{ij}-1$. Besides, the cost to update keys from the root to $L'[i]$ is $d_i$. 
More precisely, we have

\begin{align*}
    \text{cost at } L'[i] &=
    \begin{cases}
        d_i+2\sum_{j=1}^m x_{ij}-1, & \text{ if } \sum_{j=1}^m x_{ij}>0\\
        0, & \text{ if } \sum_{j=1}^m x_{ij}=0,
    \end{cases}\\
    & = d_i\times 1_{\sum_{j=1}^m x_{ij}>0}+2\sum_{j=1}^m x_{ij}\times 1_{\sum_{j=1}^m x_{ij}>0} - 1_{\sum_{j=1}^m x_{ij}>0}. 
\end{align*}

In batch insertion, with two leaf nodes being inserted new members, there is an advanced technique that we can reuse the updated key encryption keys on the overlap part of one path (the path from the root to a leaf node) to another path. It also leads to a complicated recursive relation based on the structure of the binary tree to compute the total cost. We therefore alleviate this by ignoring the overlap keys, resulting in the approximation inserting cost function as follows:
\begin{equation}
\label{ct1}
\begin{split}
   F(x) & = \sum_{i=1}^{l_1}d_i \times 1_{\sum_{j=1}^m x_{ij}>0}+2\sum_{i=1}^{l_1}\sum_{j=1}^m x_{ij}\times 1_{\sum_{j=1}^m x_{ij}>0}-\sum_{i=1}^{l_1}1_{\sum_{j=1}^m x_{ij}>0} \\
   & = \sum_{i=1}^{l_1}d_i \times 1_{\sum_{j=1}^m x_{ij}>0}+2\sum_{i=1}^{l_1}1_{\sum_{j=1}^m x_{ij}>0}\sum_{j=1}^m x_{ij}-\sum_{i=1}^{l_1}1_{\sum_{j=1}^m x_{ij}>0}\\
   & = \sum_{i=1}^{l_1}d_i \times 1_{\sum_{j=1}^m x_{ij}>0}+2\sum_{i=1}^{l_1}\sum_{j=1}^m x_{ij}-\sum_{i=1}^{l_1}1_{\sum_{j=1}^m x_{ij}>0}\\
   & = \sum_{i=1}^{l_1}(d_i-1)\times 1_{\sum_{j=1}^m x_{ij}>0}+2\sum_{i=1}^{l_1}\sum_{j=1}^m x_{ij} \\
   &= \sum_{i=1}^{l_1}{(d_i-1) \times \left \vert \sum_{j=1}^{m}{x_{ij}} \right \vert_0} + 2\sum_{i=1}^{l_1}\sum_{j=1}^m x_{ij},
\end{split}
\end{equation}
where $\vert \cdot \vert_0$ denotes the step function defined by $\vert s \vert_0 = 1$ if $s \neq 0$,  $0$ otherwise.

Moreover, we calculate the cost for a set of leaving nodes. For a given value $k\in \{1,2,...,l_2\}$, the number of new nodes $j$ being inserted at the leaf node $A[k]$ is calculated as $m_k= \sum_{j=1}^{m} y_{kj}$. It means that a subtree at the leaf node $A[k]$ includes $m_k$ leaf nodes as illustrated in Fig. \ref{fig1}. For the leaf node $A[k]$ that is chosen to append some new members (equivalent to $\sum_{j=1}^m y_{kj}>0$), the cost to build a corresponding subtree is $2\sum_{j=1}^m y_{kj}-1$. Besides, the cost to update keys from the root to $A[k]$ is $d_k-1$. 
More precisely, we have

\begin{align*}
    \text{cost at } A[k] &=
    \begin{cases}
        (d_k-1)+2\sum_{j=1}^m y_{kj}-1, & \text{ if } \sum_{j=1}^m y_{kj}>0\\
        (d_k-1), & \text{ if } \sum_{j=1}^m y_{kj}=0,
    \end{cases}\\
    & = (d_k-1) + 2\sum_{j=1}^m y_{kj}\times 1_{\sum_{j=1}^m y_{kj}>0} - 1_{\sum_{j=1}^m y_{kj}>0}. 
\end{align*}

Therefore, the approximation cost function for a set of departure nodes as follows:
\begin{equation}
\label{ct2}
\begin{split}
   F(y) & = \sum_{k=1}^{l_2}(d_k-1) +2\sum_{k=1}^{l_2}\sum_{j=1}^m y_{kj}\times 1_{\sum_{j=1}^m y_{kj}>0}-\sum_{k=1}^{l_2}1_{\sum_{j=1}^m y_{kj}>0} \\
   & = \sum_{k=1}^{l_2}(d_k-1) +2\sum_{k=1}^{l_2}1_{\sum_{j=1}^m y_{kj}>0}\sum_{j=1}^m y_{kj}-\sum_{k=1}^{l_2}1_{\sum_{j=1}^m y_{kj}>0}\\
   & = \sum_{k=1}^{l_1}(d_k-1)+2\sum_{k=1}^{l_1}\sum_{j=1}^m y_{kj}-\sum_{k=1}^{l_2}1_{\sum_{j=1}^m y_{kj}>0}\\
   &= \sum_{k=1}^{l_2}{(d_k-1) + 2\sum_{k=1}^{l_2}\sum_{j=1}^m y_{kj}-\sum_{k=1}^{l_2}\left \vert \sum_{j=1}^{m}{y_{kj}} \right \vert_0}.
\end{split}
\end{equation}
We obtain the following total cost for deleting leaving members and inserting new members 
\begin{equation}
\label{ct3}
\begin{split}
   F(x,y) & = \sum_{i=1}^{l_1}{(d_i-1) \times \left \vert \sum_{j=1}^{m}{x_{ij}} \right \vert_0} + 2\sum_{i=1}^{l_1}\sum_{j=1}^m x_{ij}\\
   &+\sum_{k=1}^{l_2}{(d_k-1) + 2\sum_{k=1}^{l_2}\sum_{j=1}^m y_{kj}-\sum_{k=1}^{l_2}\left \vert \sum_{j=1}^{m}{y_{kj}} \right \vert_0}\\
   &=\sum_{k=1}^{l_2}{(d_k-1)}+\sum_{i=1}^{l_1}{(d_i-1) \times \left \vert \sum_{j=1}^{m}{x_{ij}} \right \vert_0}-\sum_{k=1}^{l_2}\left \vert \sum_{j=1}^{m}{y_{kj}} \right \vert_0\\
   &+2\sum_{i=1}^{l_1}\sum_{j=1}^m x_{ij}+2\sum_{k=1}^{l_2}\sum_{j=1}^m y_{kj}\\
   &=\sum_{k=1}^{l_2}{(d_k-1)}+\sum_{i=1}^{l_1}{(d_i-1) \times \left \vert \sum_{j=1}^{m}{x_{ij}} \right \vert_0}-\sum_{k=1}^{l_2}\left \vert \sum_{j=1}^{m}{y_{kj}} \right \vert_0+2m.
\end{split}
\end{equation}

At the same time, we append new nodes in such a way that the balance of the tree is considered based on the given tree structure. We introduce the notion of \textit{ balance coefficient} that is defined as the difference between the index of the deepest leaf node and the index of the shallowest leaf node, $\max_{i=\overline{1,l_1},k=\overline{1,l_2}}(L'[i]\times (\sum_{j=1}^{m} x_{ij}+1),\frac{A[k]}{2}\times (\sum_{j=1}^{m} y_{kj}+1))-\min_{i=\overline{1,l_1},k=\overline{1,l_2}}(L'[i]\times (\sum_{j=1}^{m} x_{ij}+1),\frac{A[k]}{2}\times (\sum_{j=1}^{m} y_{kj}+1))$. By the above analysis, we can use the balance coefficient to control the balance of the tree: if the balance coefficient is small, the tree is encouraged to be more balanced. This coefficient is integrated into the objective function as a penalty term. It is noteworthy that the balance coefficient is much simpler than the original definition of the balance which is based on the distances from the root to the shallowest and the deepest nodes. This implies that our model will be simplified, and - consequently - is easier to handle. Finally, our optimization problem takes the following form:

\begin{align}
    \label{ct4}
    &\text{min} \sum_{i=1}^{l_1}(d_i-1)\times \biggl\vert\sum_{j=1}^{m} x_{ij}\biggr\vert_0 -\sum_{k=1}^{l_2}\biggl\vert\sum_{j=1}^{m} y_{kj}\biggr\vert_0\\
    &+\lambda\left[\text{max}_{i=\overline{1,l_1},k=\overline{1,l_2}}\left(L'[i]\times\left(\sum_{j=1}^{m} x_{ij}+1\right),\frac{A[k]}{2}\times\left(\sum_{j=1}^{m} y_{kj}+1\right)\right)\right. \notag\\
    &\left.-\text{min}_{i=\overline{1,l_1},k=\overline{1,l_2}}\left(L'[i]\times \left(\sum_{j=1}^{m} x_{ij}+1 \right),\frac{A[k]}{2}\times\left(\sum_{j=1}^{m} y_{kj}+1\right)\right)\right]\notag\\
    &\text{subject to  }\sum_{i=1}^{l_1}{x_{ij}}+\sum_{k=1}^{l_2}{y_{kj}} = 1, \forall j = \overline{1,m},\notag\\
    & \quad \quad \quad \quad \quad x_{ij} \in \{0,1\}, y_{kj}\in \{0,1\},\forall i = \overline{1,l_1},\forall k = \overline{1,l_2},\forall j = \overline{1,m} \notag,
\end{align}
where $\lambda$ is a positive parameter controlling the trade-off between the rekeying cost and the balance coefficient of the tree after insertion and deletion.

It is observed that (\ref{ct4}) is an optimization problem with binary variables and discontinuous objective. However, it can be reformulated as a combinatorial program with continuous objective by using new binary variables $u_i$ and $v_k$ as follows. Let $K = \{(x,y):\sum_{i=1}^{l_1}{x_{ij}}+\sum_{k=1}^{l_2}{y_{kj}} = 1, \forall j = \overline{1,m}, x_{ij} \in [0,1], y_{kj}\in [0,1],\forall i = \overline{1,l_1},\forall k = \overline{1,l_2},\forall j = \overline{1,m}\}$. Problem (\ref{ct4}) and the following problem have the same optimal value (refer to Proposition \ref{pro1} below)
\begin{equation}
\label{ct5}
\begin{split}
    & \text{min} \left\{\sum_{i=1}^{l_1}(d_i-1)\times u_i-\sum_{k=1}^{l_2}v_k+\lambda\left(\text{max}_{i=\overline{1,l_1},k=\overline{1,l_2}}\left(L'[i]\times(\sum_{j=1}^{m} x_{ij}+1),\frac{A[k]}{2}\times \right.\right.\right.\\
    &\left.\left.  (\sum_{j=1}^{m} y_{kj}+1)\right)-\text{min}_{i=\overline{1,l_1},k=\overline{1,l_2}}\left(L'[i]\times (\sum_{j=1}^{m} x_{ij}+1),\frac{A[k]}{2}\times(\sum_{j=1}^{m} y_{kj}+1)\right)\right): \\
    &(x,y)\in K, x \in \{0,1\}^{l_1\times m},y \in \{0,1\}^{l_2\times m},u \in \{0,1\}^{l_1},v \in \{0,1\}^{l_2}, \\
    & \left. 0\leq \sum_{j=1}^{m} x_{ij}\leq mu_i, i=\overline{1,l_1},\sum_{j=1}^{m} y_{kj}\geq v_k,k=\overline{1,l_2}\right\}
\end{split}
\end{equation}
\begin{equation}
\label{ct6}
\begin{split}
    & \Leftrightarrow\text{min} \left\{\sum_{i=1}^{l_1}(d_i-1)\times u_i-\sum_{k=1}^{l_2}v_k+\lambda\left(\text{max}_{i=\overline{1,l_1},k=\overline{1,l_2}}\left(L'[i]\times(\sum_{j=1}^{m} x_{ij}+1), \right.\right.\right.\\
    &\left.\frac{A[k]}{2}\times(\sum_{j=1}^{m} y_{kj}+1)\right)+\text{max}_{i=\overline{1,l_1},k=\overline{1,l_2}}\left(-L'[i]\times (\sum_{j=1}^{m} x_{ij}+1),-\frac{A[k]}{2}\times \right.\\
    &\left.\left.(\sum_{j=1}^{m} y_{kj}+1)\right)\right):(x,y)\in K, x \in \{0,1\}^{l_1\times m},y \in \{0,1\}^{l_2\times m},u \in \{0,1\}^{l_1}, \\
    &v \in \{0,1\}^{l_2},\left. 0\leq \sum_{j=1}^{m} x_{ij}\leq mu_i, i=\overline{1,l_1},\sum_{j=1}^{m} y_{kj}\geq v_k,k=\overline{1,l_2}\right\}.
\end{split}
\end{equation}


\begin{proposition}
\label{pro1}
For $\lambda>0$ the problems (\ref{ct4}) and (\ref{ct5}), (\ref{ct6}) are equivalent, in the sense that they have the same optimal value and $(x^*,y^*)\in K$ is a solution of (\ref{ct4}) iff there is $u^*\in \{0,1\}^{l_1},v^*\in \{0,1\}^{l_2}$ such that $(x^*,y^*,u^*,v^*)$ is a solution of (\ref{ct5}) and (\ref{ct6}).
\end{proposition}

It is observed that (\ref{ct6}) is still a very difficult optimization problem with nonsmooth objective function and binary variables. We can - however - reformulate this problem as a DC program via an exact penalty technique \citep{dinh2010efficient}. Therefore, the reformulated program can be handled by DCA which is an efficient algorithm for DC programming.

\section{Solving the batch deleting and inserting problem by DCA}\label{sec3}
\subsection{Outline of DC programming and DCA}\label{subsec3.1}

The theory of DC programming and DCA which is known to be the backbone of nonconvex programming and global optimization has been introduced first in 1985 by Pham Dinh Tao. This theory has been extensively developed from both theoretical and computational aspects since 1994 by Le Thi Hoai An and Pham Dinh Tao.

We begin by reviewing some fundamental notions and results in convex analysis. For a convex function $\theta$ defined on $\mathbb{R}^n$ and a convex set $C$, the modulus of strong convexity of $\theta$ on $C$, denoted by $\rho(\theta,C)$ or $\rho(\theta)$ if $C = \mathbb{R}^n$, is given by
\begin{equation}
\label{eqnewww}
\rho(\theta, C) = \sup\{\mu \geq 0: \theta - (\mu/2) \Vert \cdot \Vert^2 \text{ is convex on }C\}.
\end{equation}
It is said that $\theta$ is strongly convex on $C$ if $\rho(\theta,C) >0$. The subdifferential of $\theta$ at $x_0 \in \dom \theta$, denoted by $\partial \theta(x_0)$, is defined by
\begin{equation*}
\partial \theta (x_0)=\{y \in \mathbb{R}^n: \theta(x) \geq \theta(x_0) + \langle x -x_0,y \rangle, \forall x \in \mathbb{R}^n\}.
\end{equation*}
The conjugate function $\theta^{\ast}$ of $\theta \in \Gamma_0(\mathbb{R}^n)$ is defined by
$\theta^{\ast}(y)=\sup\{\langle x,y \rangle-\theta(x) : x \in \mathbb{R}^n\}.$

A function $\theta$ is said to be polyhedral convex if
\begin{align*}
    \theta(x) = \max\{\langle a_i,x \rangle + \alpha_i: \quad i=\overline{1,k}\} + \chi_C(x),
\end{align*}
where $a_i \in \mathbb{R}^n,\alpha_i \in \mathbb{R}$ for all $i=1,...,k$ and $C$ is a nonempty polyhedral convex set.

Let $\Gamma_0(\mathbb{R}^n)$ denote the convex cone of all lower semicontinuous proper convex functions on $\mathbb{R}^n$. The standard DC program takes the form
\begin{equation*}
\alpha := \inf\{f(x)= g(x) - h(x): x \in \mathbb{R}^n\} \quad (P_{dc}),
\end{equation*}
where $g,h \in \Gamma_0(\mathbb{R}^n)$. Such a function $f$ is called DC, $g-h$ is DC decomposition, while $g$ and $h$ are DC components of $f$. Note that, a DC program with closed convex constraint $x \in C$ can be equivalently written as a standard DC program in such a way that $f = (g+\chi_C) - h$, where $\chi_C$ is the indicator function of $C$.

In the DC problem $(P_{dc})$, if one of the DC components $g$ and $h$ is a polyhedral convex function, $(P_{dc})$ is called a polyhedral DC program. Polyhedral DC programs are a type of DC optimization that is used in practice and has fascinating features.

A point $x^*$ is called a critical point of $(P_{dc})$ if it verifies the generalized Karush-Kuhn-Tucker condition $\partial g(x^*) \cap \partial h(x^*) \neq \emptyset$, or equivalently $0 \in \partial g(x^*) - \partial h(x^*)$, while it is called a strongly critical point of $g-h$ if $\emptyset \neq \partial h(x^*) \subset \partial g(x^*)$.

The DCA (DC Algorithm) is an efficient technique for resolving DC programs, particularly when the DC decomposition is customised.

\textbf{The philosophy of DCA:} DCA is an iterative algorithm based on local optimality. At each iteration $k$, DCA replaces $h$ with its affine minorization by taking $y^k \in \partial h(x^k)$ and then solving the resulting convex program
\begin{equation}
(P_{k}) \quad \alpha _{k}:=\inf
\{{f}(x):={g}(x)-[h(x^k) + \langle x-x^{k},y^{k}\rangle]:x\in \mathbb{R}^n\}  \label{Pk}
\end{equation}%
whose $x^{k+1}$ is a solution, i.e., $x^{k+1}\in \partial {g}^{\ast }(y^{k})$ (see \cite{LeThi05a,PhamDinh97}).

The following standard DCA structure was inspired by the essential properties described above:

\bigskip
\noindent \textbf{Standard DCA.}

\noindent \textbf{Initialization:} Let $x^0 \in \dom \partial h$ and $k=0$.

\noindent \textbf{repeat}

Step 1: Compute the subgradient $y^{k}\in \partial h(x^{k})$.

Step 2: Solve the following convex program
\begin{equation*}
x^{k+1}\in \argmin \{g(x)-h_{k}(x):x\in \mathbb{R}^n\}.
\end{equation*}

Step 3: $k = k+1$.

\noindent \textbf{until} Stopping criterion.
\bigskip

The sequence $\{x^k\}$ generated by DCA has the following properties:
\begin{itemize}
\item[1.] The sequence $\{f(x^k)\}$ is decreasing.
\item[2.] If $f(x^{k+1})=f(x^k)$, then $x^k$ and $x^{k+1}$ are critical points of $(P_{dc})$ and DCA terminates at $k$-th iteration.
\item[3.] If the optimal value $\alpha$ of the problem $(P_{dc})$ is finite and the sequences $\{x^k\}$ and $\{y^k\}$ are bounded, then every limit point $\tilde{x}$ of $\{x^k\}$ is a critical point of $g-h.$
\item[4.] DCA has a finite convergence for polyhedral DC programs.
\end{itemize}

\subsection{A two-step DCA based algorithm for solving the problem (\ref{ct6})}\label{subsec3.2}

As the problem (\ref{ct6}) is defined on the set of leaf nodes included all leaving nodes, we first find leaf nodes in a full binary tree using the same algorithm as in \cite{le2022dc}. After finding the set of leaf nodes on the key tree, we can use the exact penalty techniques to solve the problem (\ref{ct6}). Let $p(x,y,u,v)$ be the penalty function define by
\begin{equation*}
\begin{split}
    p(x,y,u,v) &:= \sum_{i=1}^{l_1}\sum_{j=1}^{m}\text{min}\{x_{ij},1-x_{ij}\}+\sum_{k=1}^{l_2}\sum_{j=1}^{m}\text{min}\{y_{kj},1-y_{kj}\}\\
    &+\sum_{i=1}^{l_1}\text{min}\{u_i,1-u_i\}+\sum_{k=1}^{l_2}\text{min}\{v_k,1-v_k\}.
\end{split}
\end{equation*}

\noindent Notice that, there are some other possibilities to choose the penalty function. 

Then (\ref{ct6}) can be rewritten as

\begin{equation}
\label{ct7}
\begin{split}
    & \text{min} \left\{\sum_{i=1}^{l_1}(d_i-1)\times u_i-\sum_{k=1}^{l_2}v_k+\lambda\left(\text{max}_{i=\overline{1,l_1},k=\overline{1,l_2}}\left(L'[i]\times(\sum_{j=1}^{m} x_{ij}+1),\frac{A[k]}{2}\times \right.\right.\right.\\
    &\left.\left.(\sum_{j=1}^{m} y_{kj}+1)\right)+\text{max}_{i=\overline{1,l_1},k=\overline{1,l_2}}\left(-L'[i]\times (\sum_{j=1}^{m} x_{ij}+1),-\frac{A[k]}{2}\times(\sum_{j=1}^{m} y_{kj}+1)\right)\right): \\
    &(x,y)\in K, x \in [0,1]^{l_1\times m},y \in [0,1]^{l_2\times m},u \in [0,1]^{l_1},v \in [0,1]^{l_2}, 0\leq \sum_{j=1}^{m} x_{ij}\leq mu_i,\\ 
    &\left.i=\overline{1,l_1},\sum_{j=1}^{m} y_{kj}\geq v_k,k=\overline{1,l_2},p(x,y,u,v) \leq 0\right\}.
\end{split}
\end{equation}

According to \cite{dinh2010efficient}, there exists $\bar{t}$ such that the problem (\ref{ct7}) and the following  penalized problem are equivalent for all $t \geq \bar{t}$:

\begin{equation}
\label{ct8}
\begin{split}
    & \text{min} \left\{\sum_{i=1}^{l_1}(d_i-1)\times u_i-\sum_{k=1}^{l_2}v_k+\lambda\left(\text{max}_{i=\overline{1,l_1},k=\overline{1,l_2}}\left(L'[i]\times(\sum_{j=1}^{m} x_{ij}+1),\frac{A[k]}{2}\times\right.\right.\right.\\
    &\left.\left.(\sum_{j=1}^{m} y_{kj}+1)\right)+\text{max}_{i=\overline{1,l_1},k=\overline{1,l_2}}\left(-L'[i]\times (\sum_{j=1}^{m} x_{ij}+1),-\frac{A[k]}{2}\times(\sum_{j=1}^{m} y_{kj}+1)\right)\right)\\
    &+tp(x,y,u,v): (x,y)\in K, x \in [0,1]^{l_1\times m},y \in [0,1]^{l_2\times m},u \in [0,1]^{l_1},v \in [0,1]^{l_2}, \\
    &\left.0\leq \sum_{j=1}^{m} x_{ij}\leq mu_i, i=\overline{1,l_1},\sum_{j=1}^{m} y_{kj}\geq v_k,k=\overline{1,l_2}\right\}.
\end{split}
\end{equation}

Note that, in practice, the number $\bar{t}$ is generally hard to compute. Therefore, a common strategy is to use a quite large $t$ at the beginning and adaptively increase this value to encourage the solutions to be binary. Once the solution is binary, we stop increasing $t$. In particular, for our following algorithm, when a binary solution is achieved at an iteration $k$, the solution will remain binary at all latter iterations (refer to Theorem \ref{the2}).

Now let $\Delta$ be the feasible set of Problem (\ref{ct8}), i.e. $\Delta:=\{(x,y,u,v):(x,y) \in K,u \in [0,1]^{l_1},v \in [0,1]^{l_2},0 \leq \sum_{j=1}^{m} x_{ij} \leq mu_i,i=\overline{1,l_1},\sum_{j=1}^{m} y_{kj}\geq v_k,k=\overline{1,l_2}\}$ and
\begin{align*}
    \zeta(x,y,u,v) &= \lambda\left(\text{max}_{i=\overline{1,l_1},k=\overline{1,l_2}}\left(L'[i]\times(\sum_{j=1}^{m} x_{ij}+1),\frac{A[k]}{2}\times(\sum_{j=1}^{m} y_{kj}+1)\right)\right.\\
    &\left.+\text{max}_{i=\overline{1,l_1},k=\overline{1,l_2}}\left(-L'[i]\times (\sum_{j=1}^{m} x_{ij}+1),-\frac{A[k]}{2}\times(\sum_{j=1}^{m} y_{kj}+1)\right)\right).
\end{align*}
Since $p$ is concave, the following DC formulation of (\ref{ct8}) seems to be natural:
\begin{equation}
\begin{split}
    \label{ct9}
    \text{min} \{f(x,y,u,v) = &g(x,y,u,v)-h(x,y,u,v):\\
    &(x,y,u,v)\in \mathbb{R}^{l_1\times m}\times \mathbb{R}^{l_2\times m} \times\mathbb{R}^{l_1}\times\mathbb{R}^{l_2}\}
\end{split}
\end{equation}
where $g(x,y,u,v) :=\zeta(x,y,u,v)+\chi_{\Delta}(x,y,u,v)$,\\
$h(x,y,u,v) := - \sum_{i=1}^{l_1}(d_i-1)\times u_i + \sum_{k=1}^{l_2}v_k - tp(x,y,u,v),$ are clearly convex functions. In addition, the function $h$ is polyhedral convex and therefore \eqref{ct9} is a polyhedral DC program.

According to the general DCA scheme descibed above, applying DCA to (\ref{ct9}) amounts to computing two sequences $\{(x^l,y^l,u^l,v^l)\}$ and $\{(\alpha^l,\beta^l,\gamma^l,\sigma^l)\}$ in the way that $(\alpha^l,\beta^l,\gamma^l,\sigma^l)\in \partial h(x^l,y^l,u^l,v^l)$ and $(x^{l+1},y^{l+1},u^{l+1},v^{l+1})$ solves the convex program of the form $(P_{dc})$. Since $(\alpha^l,\beta^l,\gamma^l,\sigma^l)\in \partial h(x^l,y^l,u^l,v^l)$ is equivalent to
\begin{equation}
\label{ct10}
    \alpha_{ij}^l=
    \begin{cases}
        t, & \text{ if } x_{ij}^l \geq 0.5\\
        -t, & \text{ if } x_{ij}^l < 0.5,
    \end{cases}
    \forall i = 1...l_1,\forall j = 1...m,
\end{equation}
\begin{equation}
\label{ct11}
    \beta_{kj}^l=
    \begin{cases}
        t, & \text{ if } y_{kj}^l \geq 0.5\\
        -t, & \text{ if } y_{kj}^l < 0.5,
    \end{cases}
    \forall k = 1...l_2,\forall j = 1...m,
\end{equation}
\begin{equation}
    \label{ct12}
    \gamma_i^l=
    \begin{cases}
        (1-d_i)+t, & \text{ if } u_i^l \geq 0.5\\
        (1-d_i)-t, & \text{ if } u_i^l < 0.5,
    \end{cases}
    \forall i = 1...l_1,
\end{equation}
\begin{equation}
    \label{ct13}
    \sigma_k^l=
    \begin{cases}
        1+t, & \text{ if } v_k^l \geq 0.5\\
        1-t, & \text{ if } v_k^l < 0.5,
    \end{cases}
    \forall k = 1...l_2,
\end{equation}
the algorithm can be described as follow.

\bigskip

\begin{algorithm}[H]
\caption{DCAEP+ (A two-step DCA based algorithm)}

\begin{algorithmic}
\State \textbf{Step 1:} Apply the algorithm for finding leaf nodes in a full binary key tree.
\State \textbf{Step 2:}
\State \textbf{Initialization: } Let $(x^0,y^0,u^0,v^0)\in [0,1]^{l_1 \times m} \times [0,1]^{l_2 \times m} \times [0,1]^{l_1} \times [0,1]^{l_2}$ be a guess, set $l:=0,\epsilon>0,t_l>0,\theta>0.$

\Repeat
\State Compute $(\alpha^l,\beta^l,\gamma^l,\sigma^l)\in \partial h(x^l,y^l,u^l,v^l)$ via \eqref{ct10},\eqref{ct11},\eqref{ct12} and \eqref{ct13}.

\State Solve the following convex program to obtain $(x^{l+1},y^{l+1},u^{l+1},v^{l+1})$

\begin{equation}\label{ct14}
    \text{min }\{\zeta(x,y,u,v)-\langle x,\alpha^l \rangle-\langle y,\beta^l \rangle-\langle u,\gamma^l \rangle-\langle v,\sigma^l \rangle: (x,y,u,v) \in \Delta\}. 
\end{equation}
\If{$p(x^{l+1},y^{l+1},u^{l+1},v^{l+1}) > 0$} \State $t_{l+1} \leftarrow t_l+\theta.$
\EndIf

\State $l=l+1$.

\Until{$\vert f(x^{l+1},y^{l+1},u^{l+1},v^{l+1})-f(x^l,y^l,u^l,v^l)\vert \leq \epsilon(\vert f(x^l,y^l,u^l,v^l)\vert +1)$}
\end{algorithmic}
\end{algorithm}

Note that, the convex problem (\ref{ct14}) can be solved as follows. By introducing new variables $\xi = \text{max}_{i=\overline{1,l_1},k=\overline{1,l_2}}(L'[i]\times(\sum_{j=1}^{m} x_{ij}+1),\frac{A[k]}{2}\times(\sum_{j=1}^{m} y_{kj}+1)), \mu = \text{max}_{i=\overline{1,l_1},k=\overline{1,l_2}}(-L'[i]\times (\sum_{j=1}^{m} x_{ij}+1),-\frac{A[k]}{2}\times(\sum_{j=1}^{m} y_{kj}+1))$, we can reformulate (\ref{ct14}) as the following linear program which can be solved efficiently by existing optimization packages

\begin{equation}
\label{ct15}
\begin{split}
    & \text{min} \{\lambda(\xi + \mu) -\langle x,\alpha^l \rangle-\langle y,\beta^l \rangle-\langle u,\gamma^l \rangle-\langle v,\sigma^l \rangle: (x,y,u,v) \in \Delta,\\
    &L'[i]\times(\sum_{j=1}^{m} x_{ij}+1)\leq \xi,-L'[i]\times(\sum_{j=1}^{m} x_{ij}+1)\leq \mu,i=\overline{1,l_1},\\
    &\frac{A[k]}{2}\times(\sum_{j=1}^{m} y_{kj}+1)\leq \xi,-\frac{A[k]}{2}\times(\sum_{j=1}^{m} y_{kj}+1)\leq \mu,k=\overline{1,l_2}\}.
\end{split}
\end{equation}

The sequence $\{(x^l,y^l,u^l,v^l)\}$ generated by DCAEP$+$ enjoys the following property: if the penalty parameter is large enough and the solution at an iteration is binary, then the solution remains binary for the succeeding iterations. This result is interesting since DCA is a continuous approach but does have some integer property.

\begin{theorem}
\label{the2}
There exists $t^*$ such that: if $t_l > t^*$ and $(x^l,y^l,u^l,v^l)$ is binary, then $(x^s,y^s,u^s,v^s)$ remains binary for all $s \geq l$.

\end{theorem}

The convergence of DCAEP+ is stated as follows.
\begin{theorem}(Convergence properties of DCAEP+)

If at iteration $l$, $t_l>t^*$ (defined as in Theorem \ref{the2}) and $(x^l,y^l,u^l,v^l)$ is binary, 
then the DCAEP+ has a finite convergence: the sequence $\{(x^s,y^s,u^s,v^s)\}$ has a subsequence that converges to a DC critical point of $f=g-h$ after a finite number of iterations.

\end{theorem}
\section{Numerical experiments}\label{sec4}

To study the performances of our approach, we perform it on several full binary trees with different configurations. 

\subsection{Dataset}\label{subsec4.1}

We construct randomly full binary trees with the height of $8,9,10,11,12,13$ and the balance of $5,4,4,5,5,3$, respectively. First, we generate a fixed number of random integers in the appropriate value range. Secondly, we build a full binary tree with the given height based on that set of values. Therefore, the constructed trees have different configurations. A certain number of leaving members is generated randomly among a set of leaf nodes in the binary key tree.

\subsection{Comparative algorithms}\label{subsec4.2}

In \cite{le2022dc}, we proposed an optimization approach (is called as DCAEP) to an important problem in centralized dynamic group key management that consists in finding a set of leaf nodes in a binary key tree to insert new members while minimizing the insertion cost and keeping the tree as balanced as possible. The key server performs the individual deletion for the departing members and batch insertion for new members of the group. After deleting all leaving members, we use DCAEP to find a set of leaf nodes to insert new joining members.

Moreover, we will compare our optimization approach with some heuristic-based schemes that were proposed by \cite{li2001batch}, \cite{ng2007dynamic} and \cite{vijayakumar2012rotation}. We define the following variables: $D$ is the number of departing members, $J$ is the number of joining members and $H$ is the height of the key tree.

\cite{li2001batch} presented a very simple Marking algorithm that updates the key tree and generates a rekey subtree after each rekey interval time. For this algorithm, there are four cases to consider. If $J=D$, then all departing members are replaced by the joining members. If $J<D$, then they pick the $J$ shallowest leaf nodes from the departing members and replace them with the joining members. By the term “shallowest node,” they mean the leaf node of minimum height in their terminology. If $J>D$ and $D=0$, then the shallowest leaf node is selected and removed. This leaf node and the joining members form a new key tree that is then inserted at the old location of the shallowest leaf node. Next, if $J>D$ and $D > 0$, then all departing members are replaced by the joining members. The shallowest leaf node is selected from these replacements and removed from the key tree. This leaf node and the extra joining members form a new key tree that is then inserted at the old location of the removed leaf node. Last, the key server generates the necessary keys and distributes them to the members.

\cite{ng2007dynamic} devised two Merging algorithms that are suitable for batch joining events. To additionally handle
batch depart requests, they have extended these two Merging algorithms into a Batch Balanced algorithm. Their Batch Balanced Algorithm outperforms existing algorithms when the number of joining members is greater than the number of departing members and when the number of departing members is
around $N/d$ ($N$: the number of group members, $d$: the degree of the key tree) with no joining member. However, their algorithm requires the key server to multicast update messages to the members. For cases where the number of joining members and the number of departing members are comparable, their Batch Balanced Algorithm has a similar performance compared to existing work. 

\cite{vijayakumar2012rotation} presented Rotation based key tree algorithms to make the tree balanced. One limitation of their approach is that it works better for batch leave operations are more than batch join operations. When the batch join operations are more than the batch leave operations, the performance is degraded. In this case, find the insertion point to insert new nodes. It is the shallowest leaf node from the left/right subtree where the insertion of new nodes does not increase the height of the key tree. When $J=D$, replace all the leaving nodes by joining nodes. When $J<D$ and $J>0$, find all the nodes where leaving operation is going to take place. Out of the $D$ leaving nodes pick $J$ shallowest nodes in the key tree. These selected shallowest nodes are replaced by the newly joining members. Remaining $D$ nodes are simply deleted from the key tree.  

\subsection{Set up experiments and Parameters}\label{subsec4.3}

The optimization approach was implemented in the Matlab R2017b, other algorithms were implemented in the Python 3.8, and performed on a PC Intel i7-7700 CPU, 3.60GHz of 16GB RAM. CPLEX 12.6 was used for solving linear programs. We stop the DCA scheme with the tolerance $\epsilon =10^{-5}$.

Concerning the parameter $t$, as $t_0$ is hard to compute, we take a quite large value $t_0$ at the beginning and use an adaptive procedure for updating $t$ during our scheme.

\subsection{Comparative results}\label{subsec4.4}

\subsubsection{Compare with the algorithm proposed in \cite{le2022dc}}\label{subsubsec4.4.1}
\begin{table}[h!]
\begin{minipage}{\textwidth}
\begin{center}
\caption{Comparison between DCAEP and DCAEP+ in terms of rekeying
cost, balance (BL) and running time (s).}\label{tab1}
\begin{tabular}{|c|c|c|c|c|c|r|c|c|c|c|r|} 
\hline
\multirow{2}{*}{\footnotesize{$\#$}} & \multirow{2}{*}{$H$} & \multirow{2}{*}{$D$} & \multirow{2}{*}{$J$} & \multicolumn{3}{c|}{DCAEP+} & \multicolumn{5}{c|}{DCAEP}                   \\ 
&  &  & & {\footnotesize Cost} & {\footnotesize BL\footnotemark[1]}& {\footnotesize Time1} &  \thead{Insertion\\ cost} & \thead{Deletion \\cost} & \thead{Total \\cost} & {\footnotesize BL\footnotemark[1]}&{\footnotesize Time2} \\ 
\hline
 1&  8&  100&   30&  277& 3& 0.19&   118&  255&  373& 3& 0.18\\
 2&  8&  100&   50&  314& 3& 0.24&   179&  255&  434& 4& 0.22\\
 3&  8&  100&  100&  407& 4& 0.26&   293&  255&  548& 4& 0.26\\
 4&  8&  100&  200&  595& 4& 0.34&   504&  255&  759& 4& 0.41\\
 5&  8&  100&  300&  789& 4& 0.45&   703&  255&  958& 5& 0.49\\
 6&  8&  100&  400&  983& 4& 0.56&   904&  255& 1159& 5& 0.64\\
 7&  8&  100&  500& 1183& 4& 0.76&  1104&  255& 1359& 5& 0.82\\
 8&  9&  300&  200&  904& 3& 0.91&   593&  661& 1254& 4& 0.52\\
 9&  9&  300&  300& 1073& 3& 0.56&   804&  661& 1465& 4& 0.63\\
10&  9&  300&  400& 1304& 5& 1.52&  1007&  661& 1668& 5& 0.76\\
11&  9&  300&  600& 1689& 5& 1.53&  1412&  661& 2073& 5& 1.07\\
12&  9&  300&  800& 2065& 5& 1.99&  1813&  661& 2474& 5& 1.34\\
13&  9&  300& 1000& 2455& 5& 2.38&  2213&  661& 2874& 5& 1.61\\
14& 10&  700&  500& 2103& 3& 5.94&  1423& 1478& 2901& 3& 3.61\\
15& 10&  700& 1000& 3023& 6& 10.64&  2447& 1478& 3925& 6& 6.53\\
16& 10&  700& 1500& 3942& 5& 13.27&  3451& 1478& 4929& 5& 10.79\\
17& 10&  700& 2000& 4914& 5& 20.85&  4452& 1478& 5930& 5& 13.30\\
18& 10&  700& 2500& 5881& 6& 26.56&  5452& 1478& 6930& 6& 19.01\\
19& 10&  700& 3000& 6872& 6& 31.70&  6452& 1478& 7930& 6& 20.73\\
20& 11& 1000&  500& 3228& 4& 10.98&  1770& 2465& 4235& 4& 8.93\\
21& 11& 1000& 1000& 4054& 4& 18.52&  2903& 2465& 5368& 4& 16.98\\
22& 11& 1000& 1500& 4910& 6& 26.49&  3944& 2465& 6490& 6& 25.53\\
23& 11& 1000& 2000& 6271& 6& 33.88&  4953& 2465& 7418& 6& 33.33\\
24& 11& 1000& 2500& 7231& 6& 40.16&  5961& 2465& 8426& 6& 40.76\\
25& 11& 1000& 3000& 8180& 7& 49.36&  6962& 2465& 9427& 7& 48.67\\
26& 12& 3000& 2000& 8782& 4& 117.6&  5755& 6215&11970& 4& 106.7\\
27& 12& 3000& 4000&12213& 7& 228.9&  9859& 6215&16074& 7& 209.2\\
28& 12& 3000& 5000&14061& 6& 271.2& 11866& 6215&18081& 6& 265.6\\
29& 12& 3000& 6000&15867& 6& 361.2& 13868& 6215&20083& 6& 334.7\\
30& 12& 3000& 8000&19664& 7& 441.1& 17872& 6215&24087& 7& 436.8\\
31& 12& 3000&10000&23571& 7& 583.9& 21872& 6215&28087& 7& 561.9\\
32& 13& 4000& 2000&13820& 4& 366.9&  7252&10039&17291& 4& 312.9\\
33& 13& 4000& 4000&17738& 5& 650.3& 11818&10039&21857& 5& 606.4\\
34& 13& 4000& 6000&21567& 5& 1649& 16031&10039&26070& 6& 1569\\
35& 13& 4000& 8000&25373& 5& 4257& 20883&10039&30922& 5& 4135\\
36& 13& 4000& 9000&27283& 5& 6351& 22100&10039&32139& 5& 6227\\
37& 13& 4000&10000&29352& 5& 11628& 24109&10039&34148& 5& 10835\\
\hline
\end{tabular}
\end{center}
\footnotetext[1]{Balance}
\end{minipage}
\end{table}

In this problem, we consider the deletion cost and the insertion cost at the same time. Our main goal here is to find a set of leaf nodes in a binary key tree to delete leaving members and insert new members while minimizing the rekeying cost and taking the balance of the tree into account. We compare the performance of our algorithm DCAEP+ with the algorithm DCAEP proposed in \cite{le2022dc} in terms of the following three criteria: the rekeying cost (insertion cost and deletion cost), the balance of the tree after rekeying and the running time of the algorithm. The rekeying cost is the number of updated keys in actuality after deleting departure members and inserting new members. In our algorithms, it is calculated by the number of updated keys based on the approximation cost function minus the overlap keys. Moreover, the balance is measured by the difference of the distance from the root to the deepest leaf node and the shallowest one. The result is summarized in Table ~\ref{tab1}.

\textit{Comments on numerical results:}

- DCAEP$+$ always gives the lower rekeying cost than DCAEP in all cases. The cost of DCAEP$+$ on average is less than $18.6\%$, $19.6\%$, $17.9\%$, $18.0\%, 20.5\%,$ and $16.8\%$ of DCAEP in cases where the height of the tree is $8,9,10,11,12,$ and 13, respectively. 

- The balance of key tree is the same for both two algorithms.

- The running time of DCAEP+ and DCAEP is equivalent in the case of a tree whose height is $8,11,12,13$. DCAEP runs approximately $1.5$ times faster than DCAEP+ in case the height of the tree is $9$ and $10$.

\subsubsection{Compare with other existing algorithms}\label{subsubsec4.4.2}

The three mentioned heuristic-based algorithms give explicit solutions of how to delete leaving members and insert new members to the key tree, so the running time is insignificant. On the other hand, we compare the performance of DCAEP$+$ and DCAEP with other existing schemes in terms of the following two criteria: the rekeying cost (insertion cost and deletion cost) and the balance of the tree after rekeying. For each instance, we run DCAEP$+$ and DCAEP 10 times to determine the best set of parameter values.

\begin{table}[h!]
\begin{minipage}{\textwidth}
\begin{center}
\caption{Comparative results in the case of a binary tree with the height of 8, the balance of 5, and the number of leaving members of 100.}\label{tab2}
\begin{tabular}{|c|c|c|c|c|c|c|c|c|c|c|c|} 
\hline
\multirow{2}{*}{\footnotesize{$\#$}} & \multirow{2}{*}{\thead{Number \\ of new \\ members}} &   \multicolumn{2}{c|}{DCAEP$+$} & \multicolumn{2}{c|}{DCAEP} & \multicolumn{2}{c|}{Rotation} & \multicolumn{2}{c|}{Marking} & \multicolumn{2}{c|}{Merging}  \\ 
&  &   {\thead{\\Cost}} & {\footnotesize BL\footnotemark[1]} & {\thead{\\Cost}} & {\footnotesize BL\footnotemark[1]} & {\thead{\\Cost}} & {\footnotesize BL\footnotemark[1]} & {\thead{\\Cost}} & {\footnotesize BL\footnotemark[1]} & {\thead{\\Cost}} & {\footnotesize BL\footnotemark[1]} \\ 
\hline
1 & 30  & 277 & 3 & 373 & 3 & 699 & 5 & 255 &5& 258 & 5\\
2 & 50  & 314 & 3 & 434 & 4 & 699 & 5 & 255 &5& 298 & 5\\
3 & 100 & 407 & 4 & 548 & 4 & 699 & 5 & 255 &5& 398 & 5\\
4 & 150 & 496 & 4 & 656 & 4 & 827 & 1 & 354 &3& 498 & 6\\
5 & 180 & 570 & 3 & 717 & 4 & 950 & 1 & 414 &4& 558 & 8\\
6 & 200 & 595 & 4 & 759 & 4 & 997 & 1 & 454 &4& 598 & 8\\
7 & 220 & 639 & 4 & 798 & 4 & 1037& 1 & 494 &4& 638 & 8\\
8 & 250 & 693 & 4 & 859 & 5 & 1097& 1 & 554 &5& 698 & 10\\
9 & 260 & 711 & 4 & 877 & 5 & 1125& 1 & 574 &5& 718 & 10\\
10& 270 & 734 & 4 & 899 & 5 & 1145& 1 & 594 &5& 738 & 10\\
11& 280 & 746 & 4 & 918 & 5 & 1173& 1 & 614 &5& 758 & 10\\
12& 300 & 789 & 4 & 958 & 5 & 1229& 1 & 654 &5& 798 & 10\\
13& 320 & 829 & 4 & 999 & 5 & 1301& 1 & 694 &5& 838 & 10\\
14& 350 & 887 & 4 & 1059& 5 & 1369& 1 & 754 &5& 898 & 10\\
15& 370 & 921 & 4 & 1099& 5 & 1417& 1 & 794 &6& 938 & 12\\
16& 380 & 946 & 4 & 1119& 5 & 1437& 1 & 814 &6& 958 & 12\\
17& 400 & 983 & 4 & 1159& 5 & 1477& 1 & 854 &6& 998 & 12\\
18& 420 & 1019& 4 & 1199& 5 & 1517& 1 & 894 &6& 1038& 12\\
19& 450 & 1085& 4 & 1259& 5 & 1577& 1 & 954 &6& 1098& 12\\
20& 480 & 1142& 4 & 1319& 5 & 1637& 1 & 1014&6& 1158& 12\\
21& 500 & 1183& 4 & 1359& 5 & 1677& 1 & 1054&6& 1198& 12\\
\hline
\end{tabular}
\end{center}
\footnotetext [1]{Balance}
\end{minipage}
\end{table}

\begin{table}
\begin{center}
\caption{Comparative results in the case of a binary tree with the height of 9, the balance of 4, and the number of leaving members of 300.}\label{tab3}
\begin{tabular}{|c|c|c|c|c|c|c|c|c|c|c|c|} 
\hline
\multirow{2}{*}{\footnotesize{$\#$}} & \multirow{2}{*}{\thead{Number \\of new \\members}}  & \multicolumn{2}{c|}{DCAEP$+$} &\multicolumn{2}{c|}{DCAEP} & \multicolumn{2}{c|}{Rotation} & \multicolumn{2}{c|}{Marking} & \multicolumn{2}{c|}{Merging}  \\ 
&  &  {\thead{\\ Cost}} & {\footnotesize BL\footnotemark[1]} & {\thead{\\ Cost}} & {\footnotesize BL\footnotemark[1]} &{\thead{\\ Cost}} & {\footnotesize BL\footnotemark[1]} & {\thead{\\ Cost}} & {\footnotesize BL\footnotemark[1]} & {\thead{\\ Cost}} & {\footnotesize BL\footnotemark[1]} \\ 
\hline
1 & 200 & 904  & 3 & 1254 & 4 & 1959 & 4 & 661 & 4&  998& 4\\
2 & 250 & 986  & 3 & 1361 & 4 & 1959 & 4 & 661 & 4& 1098& 4\\
3 & 300 & 1073 & 3 & 1465 & 4 & 1959 & 4 & 661 & 4& 1198& 4\\
4 & 350 & 1165 & 4 & 1568 & 5 & 2090 & 1 & 760 & 5& 1298& 5\\
5 & 400 & 1304 & 5 & 1668 & 5 & 2214 & 1 & 860 & 6& 1398& 5\\
6 & 450 & 1406 & 5 & 1772 & 5 & 2354 & 1 & 960 & 7& 1498& 6\\
7 & 500 & 1500 & 5 & 1872 & 5 & 2454 & 1 & 1060& 7& 1598& 6\\
8 & 550 & 1600 & 5 & 1974 & 5 & 2554 & 1 & 1160& 7& 1698& 6\\
9 & 600 & 1689 & 5 & 2073 & 5 & 2672 & 1 & 1260& 8& 1798& 8\\
10& 650 & 1773 & 5 & 2173 & 5 & 2772 & 1 & 1360& 8& 1898& 8\\
11& 680 & 1839 & 5 & 2232 & 5 & 2832 & 1 & 1420& 8& 1958& 8\\
12& 700 & 1868 & 5 & 2273 & 5 & 2872 & 1 & 1460& 8& 1998& 8\\
13& 750 & 1955 & 6 & 2373 & 5 & 2972 & 1 & 1560& 8& 2098& 8\\
14& 780 & 2018 & 5 & 2434 & 5 & 3041 & 1 & 1620& 8& 2158& 8\\
15& 800 & 2065 & 5 & 2474 & 5 & 3081 & 1 & 1660& 8& 2198& 8\\
16& 850 & 2155 & 6 & 2573 & 5 & 3190 & 1 & 1760& 9& 2298& 10\\
17& 900 & 2251 & 5 & 2674 & 5 & 3308 & 1 & 1860& 9& 2398& 10\\
18& 950 & 2352 & 5 & 2774 & 5 & 3408 & 1 & 1960& 9& 2498& 10\\
19& 980 & 2413 & 5 & 2834 & 5 & 3468 & 1 & 2020& 9& 2558& 10\\
20& 1000& 2455 & 5 & 2874 & 5 & 3508 & 1 & 2060& 9& 2598& 10\\
\hline
\end{tabular}
\end{center}
\footnotetext[1]{Balance}
\end{table}

\begin{table}
\begin{center}
\caption{Comparative results in the case of a binary tree with the height of 10, the balance of 4, and the number of leaving members of 700.}\label{tab4}
\begin{tabular}{|c|c|c|c|c|c|c|c|c|c|c|c|} 
\hline
\multirow{2}{*}{\footnotesize{$\#$}} & \multirow{2}{*}{\thead{Number \\of new \\members}}  & \multicolumn{2}{c|}{DCAEP$+$} &\multicolumn{2}{c|}{DCAEP} & \multicolumn{2}{c|}{Rotation} & \multicolumn{2}{c|}{Marking} & \multicolumn{2}{c|}{Merging}  \\ 
&  &  {\thead{\\ Cost}} & {\footnotesize BL\footnotemark[1]} & {\thead{\\ Cost}} & {\footnotesize BL\footnotemark[1]} & {\thead{\\ Cost}} & {\footnotesize BL\footnotemark[1]} & {\thead{\\ Cost}} & {\footnotesize BL\footnotemark[1]} & {\thead{\\ Cost}} & {\footnotesize BL\footnotemark[1]} \\ 
\hline
1 & 400 & 1960 & 3 & 2675 & 3 & 4927 & 4 & 1478& 4& 2198& 4\\
2 & 500 & 2103 & 3 & 2901 & 3 & 4927 & 4 & 1478& 4& 2398& 4\\
3 & 600 & 2268 & 3 & 3109 & 3 & 4927 & 4 & 1478& 4& 2598& 4\\
4 & 700 & 2435 & 3 & 3319 & 3 & 4927 & 4 & 1478& 4& 2798& 4\\
5 & 800 & 2689 & 5 & 3520 & 5 & 3596 & 1 & 1677& 6& 2998& 5\\
6 & 900 & 2837 & 6 & 3728 & 6 & 3823 & 1 & 1877& 7& 3198& 5\\
7 & 1000& 3023 & 6 & 3925 & 6 & 4059 & 1 & 2077& 8& 3398& 6\\
8 & 1200& 3400 & 5 & 4329 & 5 & 4459 & 1 & 2477& 8& 3798& 6\\
9 & 1400& 3765 & 5 & 4729 & 5 & 4859 & 1 & 2877& 9& 4198& 8\\
10& 1500& 3942 & 5 & 4929 & 5 & 5069 & 1 & 3077& 9& 4398& 8\\
11& 1600& 4141 & 5 & 5129 & 5 & 5279 & 1 & 3277& 9& 4598& 8\\
12& 1800& 4529 & 6 & 5530 & 6 & 5679 & 1 & 3677&10& 4998& 10\\
13& 2000& 4914 & 5 & 5930 & 5 & 6079 & 1 & 4077&10& 5398& 10\\
14& 2200& 5291 & 5 & 6330 & 5 & 6479 & 1 & 4477&10& 5798& 10\\
15& 2400& 5692 & 6 & 6730 & 6 & 6879 & 1 & 4877&10& 6198& 10\\
16& 2500& 5881 & 6 & 6930 & 6 & 7079 & 1 & 5077&10& 6398& 10\\
17& 2700& 6280 & 6 & 7330 & 6 & 7479 & 1 & 5477&10& 6798& 10\\
18& 2800& 6473 & 6 & 7530 & 6 & 7679 & 1 & 5677&11& 6998& 12\\
19& 2900& 6675 & 6 & 7730 & 6 & 7879 & 1 & 5877&11& 7198& 12\\
20& 3000& 6872 & 6 & 7930 & 6 & 8079 & 1 & 6077&11& 7398& 12\\
\hline
\end{tabular}
\end{center}
\footnotetext[1]{Balance}
\end{table}

\begin{table}[h!]
\begin{center}
\caption{Comparative results in the case of a binary tree with the height of 11, the balance of 5, and the number of leaving members of 1000.}\label{tab5}
\begin{tabular}{|c|c|c|c|c|c|c|c|c|c|c|c|} 
\hline
\multirow{2}{*}{\footnotesize{$\#$}} & \multirow{2}{*}{\thead{Number \\of new \\members}}  & \multicolumn{2}{c|}{DCAEP$+$} & \multicolumn{2}{c|}{DCAEP} &\multicolumn{2}{c|}{Rotation} & \multicolumn{2}{c|}{Marking} & \multicolumn{2}{c|}{Merging}  \\ 
&  &  {\thead{ \\Cost}} & {\footnotesize BL\footnotemark[1]} & {\thead{ \\Cost}} & {\footnotesize BL\footnotemark[1]} & {\thead{ \\Cost}} & {\footnotesize BL\footnotemark[1]} & { \thead{ \\Cost}} & {\footnotesize BL\footnotemark[1]} & {\thead{ \\Cost}} & {\footnotesize BL\footnotemark[1]} \\ 
\hline
1 & 500 & 3228 & 4 & 4235 & 4 & 8865 & 5 & 2465& 5& 2998& 5\\
2 & 700 & 3550 & 4 & 4697 & 4 & 8865 & 5 & 2465& 5& 3398& 5\\
3 & 800 & 3701 & 4 & 4919 & 4 & 8865 & 5 & 2465& 5& 3598& 5\\
4 & 1000& 4054 & 4 & 5368 & 4 & 8865 & 5 & 2465& 5& 3998& 5\\
5 & 1200& 4383 & 6 & 5790 & 6 & 6203 & 1 & 2864& 7& 4398& 6\\
6 & 1500& 4910 & 6 & 6490 & 6 & 7383 & 1 & 3464& 8& 4998& 6\\
7 & 1600& 5459 & 7 & 6690 & 7 & 7616 & 1 & 3664& 9& 5198& 8\\
8 & 1800& 5844 & 7 & 7015 & 7 & 8060 & 1 & 4064& 9& 5598& 8\\
9 & 2000& 6271 & 6 & 7418 & 6 & 8592 & 1 & 4464& 9& 5998& 8\\
10& 2200& 6634 & 6 & 7824 & 6 & 9058 & 1 & 4864&10& 6398& 10\\
11& 2300& 6827 & 6 & 8021 & 6 & 9346 & 1 & 5064&10& 6598& 10\\
12& 2500& 7231 & 6 & 8426 & 6 & 9845 & 1 & 5464&10& 6998& 10\\
13& 2800& 7832 & 6 & 9026 & 6 & 10555& 1 & 6064&10& 7598& 10\\
14& 2900& 7963 & 7 & 9227 & 7 & 10832& 1 & 6264&10& 7798& 10\\
15& 3000& 8180 & 7 & 9427 & 7 & 11043& 1 & 6464&10& 7998& 10\\
16& 3200& 8590 & 6 & 9827 & 6 & 11487& 1 & 6864&11& 8398& 12\\
17& 3500& 9189 & 6 & 10427& 6 & 12230& 1 & 7464&11& 8998& 12\\
18& 3800& 9756 & 7 & 11028& 7 & 12841& 1 & 8064&11& 9598& 12\\
19& 3900& 9958 & 7 & 11228& 7 & 13041& 1 & 8264&11& 9798& 12\\
20& 4000&10167 & 7 & 11427& 7 & 13241& 1 & 8464&11& 9998& 12\\
\hline
\end{tabular}
\end{center}
\footnotetext[1]{Balance}
\end{table}

\begin{table}[h!]
\begin{center}
\caption{Comparative results in the case of a binary tree with the height of 12, the balance of 5, and the number of leaving members of 3000.}\label{tab6}
\begin{tabular}{|c|c|c|c|c|c|c|c|c|c|c|c|} 
\hline
\multirow{2}{*}{\footnotesize{$\#$}} & \multirow{2}{*}{\thead{Number \\of new \\members}}  & \multicolumn{2}{c|}{DCAEP$+$} &  \multicolumn{2}{c|}{DCAEP} & \multicolumn{2}{c|}{Rotation} &   \multicolumn{2}{c|}{Marking} & \multicolumn{2}{c|}{Merging}  \\ 
&  &  {\thead{\\ Cost}} & {\footnotesize BL\footnotemark[1]} & {\thead{\\ Cost}} & {\footnotesize BL\footnotemark[1]} & {\thead{\\ Cost}} & {\footnotesize BL\footnotemark[1]} &   {\thead{\\ Cost}} & {\footnotesize BL\footnotemark[1]} &   {\thead{\\ Cost}} & {\footnotesize BL\footnotemark[1]} \\ 
\hline
1 & 2000& 8782 & 4 & 11970 & 4 & 24258& 5 & 6215& 5& 9998& 5\\
2 & 2500& 9561 & 4 & 13023 & 4 & 24258& 5 & 6215& 5&10998& 5\\
3 & 3000&10375 & 4 & 14043 & 4 & 24258& 5 & 6215& 5&11998& 5\\
4 & 3500&11193 & 5 & 15064 & 5 & 17039& 1 & 7214& 8&12998& 6\\
5 & 4000&12213 & 7 & 16074 & 7 & 18182& 1 & 8214& 9&13998& 6\\
6 & 4200&12585 & 7 & 16472 & 7 & 18582& 1 & 8614&10&14398& 8\\
7 & 4500&13131 & 6 & 17081 & 6 & 19182& 1 & 9214&10&14998& 8\\
8 & 4800&13648 & 7 & 17684 & 7 & 19782& 1 & 9814&10&15598& 8\\
9 & 5000&14061 & 6 & 18081 & 6 & 20182& 1 &10214&10&15998& 8\\
10& 5500&14929 & 7 & 19082 & 7 & 21182& 1 &11214&11&16998& 10\\
11& 5800&15510 & 6 & 19684 & 6 & 21818& 1 &11814&11&17598& 10\\
12& 6000&15867 & 6 & 20083 & 6 & 22230& 1 &12214&11&17998& 10\\
13& 6500&16803 & 7 & 21084 & 7 & 23230& 1 &13214&11&18998& 10\\
14& 7000&17771 & 6 & 22087 & 6 & 24242& 1 &14214&11&19998& 10\\
15& 7500&18760 & 6 & 23086 & 6 & 25254& 1 &15214&12&20998& 12\\
16& 8000&19664 & 7 & 24087 & 7 & 26266& 1 &16214&12&21998& 12\\
17& 8500&20650 & 6 & 25087 & 6 & 27266& 1 &17214&12&22998& 12\\
18& 9000&21595 & 7 & 26087 & 7 & 28266& 1 &18214&12&23998& 12\\
19& 9500&22575 & 7 & 27086 & 7 & 29266& 1 &19214&12&24998& 12\\
20&10000&23571 & 7 & 28087 & 7 & 30266& 1 &20214&12&25998& 12\\
\hline
\end{tabular}
\end{center}
\footnotetext[1]{Balance}
\end{table}

The specific test scenario is executed as follows: Deleting a certain number of leaving members and inserting several numbers of new members in full binary trees that have different heights according to the chosen algorithm. We compare the performance of DCAEP$+$ and DCAEP with three existing schemes: the Rotation algorithm, the Marking algorithm and the Merging algorithm. The rekeying cost and balance of the original trees that have the height of $8,9,10,11$ and 12, and the balance of $5,4,4,5,5$ are summarized in Tables~\ref{tab2}--\ref{tab6}, respectively.

\textit{Comments on numerical results:}

- In terms of updating key cost, DCAEP$+$ always results in a lower rekeying cost than the Rotation algorithm in all cases. The cost of DCAEP$+$ on average is less than $36.3\%$, $36.4\%$, $25.4\%$, $32.1\%,$ and $32.6\%$ of the Rotation algorithm corresponding to the height of the tree at $8,9,10,11,$ and 12, respectively. DCAEP also has a lower cost than the Rotation algorithm, namely $22.8\%$, $21.9\%$, $8.7\%$, $19.6\%,$ and $15.1\%$ less for trees with heights of $8,9,10,11,$ and 12, respectively. For the case where the number of new members greater than the number of leaving members $(J>D)$, DCAEP$+$ has the same cost with the Merging algorithm. While the Marking algorithm gives the lowest updating key cost, it cannot maintain the tree balance after deletion and insertion. The rekeying cost of DCAEP is slightly higher than that of the Marking and Merging algorithms.

- As for the balance of the key tree after deletion and insertion, DCAEP$+$ always keeps this factor lower than the Marking algorithm. In particular, the balance of our approach is less than $24.3\%$, $34.5\%$, $39.5\%$, $32.8\%$ and $44.8\%$ of the Marking algorithm where the tree has a height of $8,9,10,11$ and 12, respectively. DCAEP$+$ is more efficient in terms of balance than the Merging algorithm in the case of a binary tree with the height of 8. In these instances, DCAEP$+$ has an average balance of $5$ less than $59.3\%$ of the Merging algorithm's. In the case of a binary tree with the height of $9,10,11,12$, where the number of new members is greater than or equal to $450,1200,1600,4200$ and greater than the number of departure members (on $15,13,14,15/20$ instances), our algorithm is more efficient in terms of balance than the Merging algorithm. In these instances, DCAEP$+$ has an average balance of $5.1, 5.5, 6.5,6.5$, less than $37.9\%, 42.9\%, 36.8\%,36.4\%$ of the Merging algorithm's where the tree has a height of $9,10,11$ and 12, respectively. The Marking and Merging algorithms keep the balance of the tree after rekeying as same as the original tree when the number of joining members is less than or equal to the number of leaving members. Moreover, DCAEP+ and DCAEP keep the same balance of the key tree after rekeying.

Overall, our proposed algorithm takes into account simultaneously both objectives: the rekeying cost and the balance of the tree. It is efficient in both scenarios when the number of leaving members exceeds the number of new members and vice versa for key trees that are more unbalanced.

\section{Conclusion}\label{sec5}

In this paper, we have proposed an optimization approach to the problem of batch deletion and insertion members in the LKH structure. Our primary objective is to simultaneously minimize the cost of removing and adding nodes while maintaining a balanced tree. The suggested optimization problem has binary variables and an objective function that is discontinuous. It is first equivalently formulated to eliminate the objective's step functions. By using recent results on exact penalty techniques in DC programming, the latter problem can be reformulated as a DC program. We then designed an efficient DCA for solving this problem. Numerical experiments have been done to demonstrate the validity of our suggested model and its corresponding DCA. In comparison to existing approaches, it has been shown that our approach achieves a better trade-off between two considered criteria. To conclude, our method allows for a good compromise between the rekeying cost, which is the cost of deleting and inserting members, and the tree's balance after rekeying.

\backmatter


\bibliography{sn-bibliography}


\end{document}